\documentstyle[12pt,amssymb,graphicx]{article}
\def\lk{{\mathrm lk}}
\date{} 
\renewcommand{\thefootnote}{}
\begin{document}

\def\thebibliography#1{\begin{center}{\normalsize\bf References}\end{center}
 \list{[{\bf \arabic{enumi}}]}{\settowidth\labelwidth{[#1]}
 \leftmargin\labelwidth
 \advance\leftmargin\labelsep
 \usecounter{enumi}}
 \def\newblock{\hskip .11em plus .33em minus .07em}
 \sloppy\clubpenalty4000\widowpenalty4000
 \sfcode`\.=1000\relax}

\title{
{\large LINKING NUMBERS IN RATIONAL HOMOLOGY $3$-SPHERES, 
CYCLIC BRANCHED COVERS AND INFINITE CYCLIC COVERS}}

\author{{\normalsize J\'{O}ZEF H. PRZYTYCKI}\\
{\small Department of Mathematics, The George Washington University}\\[-1mm]
{\small Washington, DC 20052, USA}\\[-1mm]
{\small e-mail: przytyck@research.circ.gwu.edu}\\[3mm]
{\normalsize AKIRA YASUHARA }\\
{\small Department of Mathematics, Tokyo Gakugei University}\\[-1mm]
{\small Nukuikita 4-1-1, Koganei, Tokyo 184-8501, Japan}\\[-1mm]
{\small e-mail: yasuhara@u-gakugei.ac.jp}\\
}

\maketitle

\baselineskip=16pt
\vspace*{-5mm}  
{\small 
\begin{quote}
\begin{center}A{\sc bstract}\end{center} 
We study the linking numbers in 
a rational homology $3$-sphere and 
in the infinite cyclic cover of the complement of a knot. 
They take values in $\Bbb Q$ and in  ${Q}({\Bbb Z}[t,t^{-1}])$ 
respectively, where ${Q}({\Bbb Z}[t,t^{-1}])$ denotes the quotient 
field of ${\Bbb Z}[t,t^{-1}]$. 
It is known that the modulo-$\Bbb Z$ linking number in 
the rational homology $3$-sphere is determined by the linking matrix 
of the framed link 
and that the modulo-${\Bbb Z}[t,t^{-1}]$ 
linking number in the infinite cyclic cover of the complement of a knot
is determined by the Seifert matrix of the knot. 
We eliminate \lq modulo $\Bbb Z$' and 
\lq modulo ${\Bbb Z}[t,t^{-1}]$'. 
When the finite cyclic cover  
of the $3$-sphere  branched over a knot is a rational homology $3$-sphere, 
the linking number of a pair in the preimage of 
a link in the $3$-sphere  
is determined by the Goeritz/Seifert matrix of the knot.
\end{quote}}


\renewcommand{\thefootnote}{*}

\baselineskip=20pt

\newpage
\bigskip
{\bf Introduction}

\medskip

Let $K\cup K_1\cup\cdots\cup K_m$ $(m\geq 1)$ be an oriented 
$(m+1)$-component link in the three sphere $S^3$. 
If the linking number $\lk(K,K_i)$ is even 
for any $i(=1,...,m)$, then there is an unoriented, 
possibly nonorientable surface $F$ bounded 
by $K$ disjoint from $K_1\cup\cdots\cup K_m$. 
Let $G_{\alpha}$ be the {\em Goeritz matrix} \cite{Goe}, \cite{G-L}
with respect to a basis
$\alpha=(a_1,...,a_n)$ of $H_1(F)$, i.e., 
the $(i,j)$-entry of $G_{\alpha}$ is 
equal to $\lk(a_i,\tau a_j)$, where $\tau a_j$ is a 1-cycle 
in $S^3-F$ obtained by pushing off $2a_j$ into both normal 
directions.\footnote{$2a_j$ can be thought as the double 
cover of $a_j$ lying in the boundary of the regular neighborhood of $F$.} 
Let $V_{\alpha}(K_i)=(\lk(K_i,a_1),...,\lk(K_i,a_n))$. 
For $i,j$ $(1\leq i,j\leq m,$ possibly $i=j$) 
we define 
\[\lambda_F(K_i,K_j)=
V_{\alpha}(K_i)G_{\alpha}^{-1}V_{\alpha}(K_j)^T,\]
and $\lambda_F(K_i,K_j)=0$ for a $2$-disk $F$. 
Note that $\lambda_F(K_i,K_j)=\lambda_F(K_j,K_i)$.
The number $\lambda_F(K_i,K_j)$ is independent of the choice of a basis 
and $S^*$-equivalence calss of $F$ in $S^3-(K_i\cup K_j)$ 
(Proposition 2.1), and if $i=j$ it is an invariant of links (Corollary 2.4). 

If $\lk(K,K_i)=0$ for any $i(=1,...,m)$, then there is 
a Seifert surface $F$ of $K$ with $F\cap (K_1\cup\cdots\cup K_m) 
=\emptyset$. Let $M_{\alpha}$ be the {\em Seifert matrix} with respect to 
a basis $\alpha=(a_1,...,a_n)$, 
i.e., $m_{ij}=\lk(a^+_i,a_j)(=
\lk(a_i,a^-_j))$, 
where $a_i^{\pm}$ means a curve that is obtained by pushing 
off into $\pm$-direction.  
Let $G_{\alpha,\omega}$ be a Hernitean matrix 
$(1-\overline{\omega})M_{\alpha}+(1-\omega)M_{\alpha}^T$, 
where $\omega(\neq1)$ is a root of unity different from 
a root of the Alexander polynomial of $K$. 
Since 
$G_{\alpha,\omega}=(\omega-1)(\overline{\omega}M_{\alpha}-M_{\alpha}^T$), 
$G_{\alpha,\omega}$ is nonsingular. 
For $i,j$ $(1\leq i,j\leq m,$ possibly $i=j)$ we define 
\[\lambda_F(K_i,K_j;\omega)=
V_{\alpha}(K_i)G_{\alpha,\omega}^{-1}V_{\alpha}(K_j)^T,\]
and $\lambda_F(K_i,K_j;\omega)=0$ for a $2$-disk $F$. 
Let $G_{\alpha}(t)=tM_{\alpha}-M_{\alpha}^T$. Then we define 
\[\lambda_F(K_i,K_j)(t)=
V_{\alpha}(K_i)G_{\alpha}(t)^{-1}V_{\alpha}(K_j)^T,\]
and $\lambda_F(K_i,K_j)(t)=0$ for a $2$-disk $F$.
Let $M_{p,\alpha}$ be a $(p-1)n\times (p-1)n$ matrix defined by 
\[M_{p,\alpha}=\left(
\begin{array}{ccccc}
M_{\alpha}+M_{\alpha}^T & -M_{\alpha}^T & O & \cdots & O \\
-M_{\alpha} & M_{\alpha}+M_{\alpha}^T & \ddots & \ddots & \vdots \\
O & \ddots & \ddots & \ddots & O \\
\vdots & \ddots & \ddots & M_{\alpha}+M_{\alpha}^T & -M_{\alpha}^T \\
O & \cdots & O & -M_{\alpha} & M_{\alpha}+M_{\alpha}^T
\end{array}\right),\]
where $O$ is the $n\times n$ zero matrix. 
Note that $M_{p,\alpha}$ is a presentation matrix of the first homology 
group of the $p$-fold cyclic cover of $S^3$ branched over $K$ 
\cite{Kau}.
Let $V_{p,\alpha}^k(K_i)=({\bf 0}_{(k-1)n},V_{\alpha}(K_i),
{\bf 0}_{(p-k-1)n})$, 
where ${\bf 0}_l$ is the $1\times l$ zero vector. 
When $M_{p,\alpha}$ is nonsingular, i.e., 
the $p$-fold cyclic cover of $S^3$ branched over $K$ is 
a rational homology $3$-sphere,  
we define
\[\lambda_F^{(k,l)}(K_i,K_j)=
V_{p,\alpha}^k(K_i)M_{p,\alpha}^{-1}V_{p,\alpha}^l(K_j)^T,\]
and $\lambda_F^{(k,l)}(K_i,K_j)=0$ for a $2$-disk $F$.
Note that $\lambda_F(K_i,K_j;\omega)=\lambda_F(K_j,K_i;\omega)$, 
$\lambda_F^{(k,l)}(K_i,K_j)=\lambda_F^{(k,l)}(K_j,K_i)$, and   
$\lambda_F(K_i,K_j)(t)$  is equal to $\lambda_F(K_j,K_i)(t^{-1})$ 
up to multiplication by  a unit of ${\Bbb Z}[t,t^{-1}]$. 
We shall show that $\lambda_F(K_i,K_j;\omega)$, 
$\lambda_F(K_i,K_j)(t)$ and $\lambda^{(k,l)}_F(K_i,K_j)$
are independent of the choice of basis and 
$S$-equivalence class of $F$ in $S^3-(K_i\cup K_j)$ 
(Proposition 3.1), and 
if $i=j$ then $\lambda_F(K_i,K_i;\omega)$ and 
$\lambda_F(K_i,K_i)(t)$ are invariants for links (Corollary 4.2). 
The definitions of $\lambda_F(K_i,K_j)$ etc. were given by Y.W. Lee 
\cite{Lee1}, \cite{Lee2}. But his definitions require some additional 
condition. We make his definitions more general. 

Let $M$ be a rational homology $3$-sphere and $K_1\cup K_2$ a 
$2$-component oriented link in $M$. Then there is a $2$-chain $F$ in 
$M$ such that $F$ bounds $cK_1$, where $cK_1$ is a disjoint 
union of $c$ copies of $K_1$ in a small 
neighborhood of $K_1$.  We define 
\[\lk_M(K_1,K_2)=\frac{F\cdot K_2}{c}\in {\Bbb Q},\]
where $F\cdot K_2$ is the intersection number of $F$ and $K_2$ \cite{S-T}. 
It is known that this linking number is well-defined and 
$\lk_M(K_1,K_2)=\lk_M(K_2,K_1)$. 
Note that $\lk_{S^3}$ is as same as the linking number 
$\lk$ in the usual sense. 

Let $S$ be a $3$-manifold with the boundary composed of some tori. 
Let $S_{\mu}$ and $S_{\delta}$ be $3$-manifold 
obtained from $S$ by Dehn fillings with respect to systems of curves
$\mu$ and $\delta$ on $\partial S$ respectively. 
Suppose that both $S_{\mu}$ and $S_{\delta}$ are 
rational homology 3-spheres. 
In Section 1, we show that the difference of the linking 
number $\lk_{S_{\delta}}-\lk_{S_{\mu}}$ is determined by 
a matrix obtained from  $\mu$ and $\delta$  (Theorem 1.1). 
It generalizes a result of J. Hoste \cite{Hos} proved for  
{\em integral} homology $3$-spheres.
As a corollary, for a rational homology $3$-sphere $M$ obtained 
by Dehn surgery along a framed link in $S^3$, we obtain that 
the linking number $\lk_M$ is determined by the linking matrix 
of the framed link (Corollary 1.2). 
It is known that the linking number modulo $\Bbb Z$ is obtained via 
the matrix; see \cite{G-L} for example. 
Our results does not require \lq modulo $\Bbb Z$'.

In sectins 2 and 3 we show that, 
for a 3-component link $K\cup K_1\cup K_2$ with 
$\lk(K,K_i)$ even (resp. $=0$),  
$\lk_{X_{2}}(K_{ik},K_{jl})$ (resp. $\lk_{X_{p}}(K_{ik},K_{jl})$) 
is determined by $\lambda_F(K_i,K_j)$ (resp. $\lambda_F^{(k,l)}(K_i,K_j)$) 
(Theorems 2.3 and 3.2), 
where $X_p$ is the $p$-fold cyclic cover of $S^3$ 
branched over $K$,  and $K_{ik}(\subset X_p)$ is a component of 
the preimage of $K_i$.  

Let $X_{\infty}$ be the infinite cyclic cover of the complement 
of a knot and $\tau$ a covering translation that 
shifts $X_{\infty}$ along the positive direction 
with respect to the knot. 
Let $K_1\cup K_2$ be a 
$2$-component oriented link in $X_{\infty}$ with 
$\tau^i K_1\cap K_2=\emptyset$ for any $i\in{\Bbb Z}$. 
Note that there is a $2$-chain $F$ in $X_{\infty}$ such that 
\[\partial F=\bigcup_{k\in{\Bbb Z}}c_{k}\tau^k K_1,\]
where $c_{k}$'s are integers. Then we define 
\[\tilde{\lk}_{X_{\infty}}(K_1,K_2)=
\frac{\sum_{h\in {\Bbb Z}}t^h(F\cdot \tau^h K_2)}{\sum_{k\in{\Bbb Z}}c_{k}t^k}
\in{Q}({\Bbb Z}[t,t^{-1}]).
\]
Since $H_2(X_{\infty};{\Bbb Z})\cong 0$ \cite{B-Z}, 
this is well-defined. 
We do not need to treat this linking pairing  
up to modulo ${\Bbb Z}[t,t^{-1}]$.
Note that $\tilde{\lk}_{X_{\infty}}(\tau K_i,K_j)=
\tilde{\lk}_{X_{\infty}}(K_i,\tau^{-1} K_j)=
t\tilde{\lk}_{X_{\infty}}(K_i,K_j)$. 
For a parallel copy $K'_i$ of $K_i$ with $\lk(K_i,K'_i)=0$, 
the linking pairing $\tilde{\lk}_{X_{\infty}}(K_i,K'_i)$ is called  
{\em Kojima-Yamazaki's $\eta$-function} $\eta(K,K_i;t)$ \cite{K-Y}.

In Section 4 we show that, for a 3-component link $K\cup K_1\cup K_2$ 
with $\lk(K,K_i)=0$, 
$\tilde{\lk}_{X_{\infty}}({K_{ik}},{K_{jl}})$ is determined by 
$\lambda_F(K_i,K_j)(t)$ (Therem 4.1). This means 
the linking pairing is obtained via the Seifert 
matrix of $K$. It is known that the linking pairing 
modulo ${\Bbb Z}[t,t^{-1}]$ is determined by the matrix  
\cite{Kea}, \cite{Lev}, \cite{Tro2}. 
Our result does not require \lq modulo ${\Bbb Z}[t,t^{-1}]$'.  
As a corollary we 
have that $(1-t)\lambda_F(K_i,K_i)(t)$ is equal to 
Kojima-Yamasaki's  $\eta$-function $\eta(K,K_i;t)$ and that 
$(1-t)\lambda_F(K_1,K_2)(t)+\lk(K_1,K_2)$ is a topological 
concordance invariant of 
$K\cup K_1\cup K_2$ up to multiplication by $t^{\pm n}$. 

\bigskip
{\bf 1. Rational homology $3$-sphere}

\medskip
Let $S$ be a 3-manifold with a boundary composed of $n$ tori,
$T^2_1,T^2_2,...,T_n^2$.
Suppose that $\mu =(\mu_1,\mu_2,...,\mu_n)^T$ and 
$\delta =(\delta_1,\delta_2,...,\delta_n)^T$ are two systems of curves
(written as columns),
$\mu_i,\delta_i \subset T^2_i$, such that the intersection number 
$\mu_i \cdot \delta_i = q_i \neq 0$. Furthermore we suppose that
{\em$\mu$ and $\delta$ represent two bases of $H_1(S;{\Bbb Q})$}.
This condition can be restated as: 
Dehn fillings of $S$ with respect to $\mu$ and $\delta$ 
give rational homology spheres $S_{\mu}$ and $S_{\delta}$ respectively.
Since $\mu$ and $\delta$ represent two bases 
$[\mu] =([\mu_1],[\mu_2],...,[\mu_n])^T$ and 
$[\delta] =([\delta_1],[\delta_2],...,[\delta_n])^T$of $H_1(S;{\Bbb Q})$, 
there is an $n\times n$-matrix $B=(b_{ij})$ changing the basis, which is 
an invertible matrix with rational coefficients such that 
$\delta_i=\sum_{j=1}^n b_{ij}\mu_j$ or shortly 
$[\delta] = B[\mu]$ (and $[\mu] = B^{-1}[\delta]$).
Let $J_i$ (resp. $\widehat{J_i}$) be the core of a solid torus 
attached to $T_i^2$ in $S_{\mu}$ (resp. $S_{\delta}$). 
Let $G=(g_{ij})$ be an $n\times n$-matrix with 
$g_{ij}=\lk_{S_{\mu}}(J_i,J_j)$ for $i\neq j$  
and $g_{i,i}= {b_{ii}}/{q_i}$. Note that  
${b_{ii}}/{q_i}$ is a Dehn surgery coefficient used to 
change $S_{\mu}$ to $S_{\delta}$.
In particular $[\delta_i - b_{ii}\mu_i]$ is zero in 
$H_1(S_{\mu} - J_i;{\Bbb Q})$. We call $G=(g_{ij})$ a 
{\em surgery-linking matrix} from $S_{\mu}$ to $S_{\delta}$. 
We can consider the surgery-linking matrix $H=(h_{ij})$ from
$S_{\delta}$ to $S_{\mu}$ in an analogous manner, i.e., 
$h_{ij}=\lk_{S_{\delta}}(\widehat{J_i},\widehat{J_j})$ for $i\neq j$
and $h_{ii}=\overline{b_{ii}}/{(-q_i)}$, where 
$\overline{b_{ij}}$ is the $(i,j)$-entry of $B^{-1}$. 
Note that $q_i=\mu_i\cdot \delta_i=-\delta_i\cdot\mu_i$.
Let $Q$ be a diagonal matrix with $q_{ii}=q_i$. Then we have 
the following theorem. 

\medskip
{\bf Theorem 1.1.} {\em 
\begin{enumerate}
\item[{\rm (1)}] $B=QG$ and $B^{-1}=-QH$.
\item[{\rm (2)}] For a two component oriented link $K_1\cup K_2$ in $S$, 
\[\begin{array}{l}
\lk_{S_{\delta}}(K_1,K_2) - \lk_{S_{\mu}}(K_1,K_2)\\
\hspace*{3cm}=-(\lk_{S_{\mu}}(K_1, J_1),...,\lk_{S_{\mu}}(K_1,J_n))G^{-1}
(\lk_{S_{\mu}}(K_2,J_1),...,\lk_{S_{\mu}}(K_2, J_n))^T.
\end{array}
\]\end{enumerate}
}

\medskip 
In Theorem 1.1(2), the case that both 
$S_{\mu}$ and $S_{\delta}$ are {\em integral} 
homology $3$-spheres was shown by J. Hoste \cite{Hos}. 

Before proving Theorem 1.1, we formulate a useful corollary. 
Let $J_1\cup\cdots \cup J_n$ be an $n$-component oriented link in $S^3$.
We say that $J_1\cup\cdots \cup J_n$ is a {\em $($rational$)$ framed link} 
if every component $J_i$ is equipped with a rational number  $p_i/q_i$ 
with $p_i,q_i\in{\Bbb Z}$. 
Let $N_i$ be a small neighborhood of $J_i$ in $S^3$ such that $N_i\cap N_j=
\emptyset$ for $i\neq j$. Let $m_i$ be a meridian of $N_i$ 
with $\lk(m_i,J_i)=1$ and $l_i$ a longitude that is null-homologous in 
$S^3-J_i$. 
Then we obtaine a new $3$-manifold $M$ in the following way: 
Remove the interiors of the tori $N_1,...,N_n$ from $S^3$, attach 
$2$-handles $D^2_1\times [0,1],...,D^2_n\times [0,1]$  
so that $[\partial D_i]=p_i[m_i]+q_i[l_i]\in H_1(\partial N_i)$ $(i=1,...,n)$, 
and cap off it with $3$-balls. We say that $M$ is 
{\em obtained by Dehn surgery along the $($rational$)$ framed link} 
$J_1\cup\cdots \cup J_n$. 
Let $ {G}=(g_{ij})$ be the {\em linking matrix} of the framed link, i.e., 
$g_{ij}=\lk_{S^3}(J_i,J_j)$ if $i\neq j$ and $g_{ii}=p_i/q_i$. 
Since $G$ is a surgery-linking matrix from 
$S^3$ to $M$, by Theorem 1.1(2), we have the following corollary. 

\medskip
{\bf Corollary 1.2.} {\em 
Let $M$ be a rational homology $3$-sphere obtained by Dehn surgery 
along a rational framed, oriented link 
$J_1\cup\cdots \cup J_n$ in $S^3$. 
Let ${G}$ be the linking matrix of the framed link. 
Then for a $2$-component oriented link $K_1\cup K_2$ in the 
complement of the framed link, 
\[\begin{array}{l}
\lk_M(K_1,K_2)-\lk_{S^3}(K_1,K_2)\\
\hspace*{3cm}=-(\lk_{S^3}(K_1,J_1),...,
\lk_{S^3}(K_1,J_n)){G}^{-1}
(\lk_{S^3}(K_2,J_1),...,\lk_{S^3}(K_2,J_n))^T.\ \Box
\end{array}\]}

\medskip
{\bf Proof of Theorem 1.1.} 
(1) By the definitions of $G$ and $H$, we have $b_{ii}=q_ig_{ii}$ and 
$\overline{b_{ii}}=-q_ih_{ii}$. 
We may assume $i\neq j$. 
Since $[\delta]=B[\mu]$, each $d\delta_i$ is homologous to 
$d\sum_{k=1}^n b_{ik}\mu_k$ in $S$ for some integer $d$. 
This implies 
$\lk_{S_{\mu}}(\delta_i,J_j)=\lk_{S_{\mu}}(\sum_{k=1}^n b_{ik}\mu_k,J_j)
=b_{ij}$. 
Meanwhile $\delta_i$ is homologous to $q_iJ_i$ in the solid torus attached 
$T^2_i$ since $\mu_i\cdot \delta_i=q_i$. Therefore 
$\lk_{S_{\mu}}(\delta_i,J_j)=\lk_{S_{\mu}}(q_iJ_i,J_j)=q_ig_{ij}$.
Notice that $\delta_i \cdot \mu_i =-\mu_i \cdot \delta_i = -q_i$. 
By the same arguments as in above, we have $\overline{b_{ij}}=
\lk_{S_{\delta}}(\sum_{k=1}^n \overline{b_{ik}}\delta_k,\widehat{J_j})
=\lk_{S_{\delta}}(\mu_i,\widehat{J_j})
=\lk_{S_{\delta}}(-q_i\widehat{J_i},\widehat{J_j})=-q_ih_{ij}$. 

(2) Since $dK_k$ is homologous to 
$d\sum_{i=1}^n \lk_{S_{\mu}}(K_k,{J_i})\mu_i$ in $S$ for some integer $d$, 
there is a 2-chain $F_k$ in $S$ that realizes the homologous above. 
This implies that 
\[\begin{array}{rcl}
\displaystyle
\lk_{S_{\delta}}(K_1,K_2)-\frac{F_1\cdot K_2}{d}&=&\displaystyle
\lk_{S_{\delta}}\left(\sum_{i=1}^n \lk_{S_{\mu}}(K_1,
{J_i})\mu_i,K_2\right)\\
&=&\displaystyle
\lk_{S_{\delta}}\left(\sum_{i=1}^n \lk_{S_{\mu}}(K_1,{J_i})\mu_i, 
\sum_{j=1}^n \lk_{S_{\mu}}(K_2,{J_j})\mu_j\right)\\
&=&\displaystyle
\sum_{i=1}^n\sum_{j=1}^n \lk_{S_{\mu}}(K_1,{J_i})\lk_{S_{\mu}}(K_2,{J_j})
\lk_{S_{\delta}}(\mu_i,\mu_j)\end{array}\]
 and 
\[\begin{array}{rcl}
\displaystyle
\lk_{S_{\mu}}(K_1,K_2)-\frac{F_1\cdot K_2}{d}&=&\displaystyle
\lk_{S_{\mu}}\left(\sum_{i=1}^n \lk_{S_{\mu}}(K_1,{J_i})\mu_i, 
K_2\right)\\
&=&\displaystyle
\lk_{S_{\mu}}\left(\sum_{i=1}^n \lk_{S_{\mu}}(K_1,{J_i})\mu_i, 
\sum_{j=1}^n \lk_{S_{\mu}}(K_2,{J_j})\mu_j\right)\\
&=&0,\end{array}\]
where $\lk_{S_{\delta}}(\mu_i,\mu_i)$ and $\lk_{S_{\mu}}(\mu_i,\mu_i)$ 
mean the linking numbers of $\mu_i$ and a parallel copy of 
$\mu_i$ in $T_i$. Hence we have 
\[\lk_{S_{\delta}}(K_1,K_2) - \lk_{S_{\mu}}(K_1,K_2) =
\sum_{i=1}^n\sum_{j=1}^n \lk_{S_{\mu}}(K_1,{J_i})
\lk_{S_{\mu}}(K_2,{J_j})\lk_{S_{\delta}}(\mu_i,\mu_j).\] 
Since $\mu_i$ is homologous to $-q_i\widehat{J_i}$ 
in the solid torus attached to $T^2_i$ and since 
$d\mu_i$ is homologous to 
$d\sum_{k=1}^n \overline{b_{ik}}\delta_k$ in $S$ for some integer $d$, 
$\lk_{S_{\delta}}(\mu_i,\mu_j)=
\lk_{S_{\delta}}(-q_i\widehat{J_i},-q_j\widehat{J_j})=q_iq_jh_{ij}$ 
for $i\neq j$ and 
$\lk_{S_{\delta}}(\mu_i,\mu_i)=\lk_{S_{\delta}}(-q_i\widehat{J_i}, 
\sum_{k=1}^n \overline{b_{ik}}\delta_k)=-q_i\overline{b_{ii}}=
q_iq_ih_{ii}$. So we have 
$\lk_{S_{\delta}}(\mu_i,\mu_j)=q_iq_jh_{ij}$ for any $i,j$. 
Theorem 1.1(1) completes the proof. 
$\Box$

\medskip
{\bf Remark 1.3.}
In Theorem 1.1(1), the assumption that both $[\mu]$ and $[\delta]$ 
are bases of $H_1(S;{\Bbb Q})$ is not necessarily needed. 
We can obtain the same result if $[\mu]$ and $[\delta]$ are bases of 
the same subspace of $H_1(S;{\Bbb Q})$.  

\medskip
{\bf Remark 1.4.}
(a) Let $K_1\cup K_2$ be a $2$-component oriented link in 
an oriented manifold $M$ 
each of which component represents an element 
in Tor$H_1(M)$. For a $2$-chain $F$ in 
$M$ with $\partial F=cK_1$, we define 
\[\lk_M(K_1,K_2)=\frac{F\cdot K_2}{c}\in {\Bbb Q}.\]
Since $[K_2]$ is in Tor$H_1(M)$, 
$(c'F\cup (-cF'))\cdot K_2= 0$ 
for any 2-chain $F'$ with $\partial F'=c'K_1$. 
This implies that $\lk_M$ is well-defined. 

(b) Let $M$ be an oriented $3$-manifold. 
We define a function mul$:H_1(M)\rightarrow {\Bbb Z}$ as 
follows: For an element $a\in H_1(M)$, let mul$(a)$ be the 
greatest common divisor of the integers in 
$\{a\cdot F|F\mbox{ is a $2$-cycle in }M\}$. We put mul$(a)$=0 
if $a\cdot F=0$ for any $F$. 
Set $T(H_1(M))=\{a\in H_1(M)|{\mathrm mul}(a)=0\}$. 
Note that Tor$(H_1(M))\subset T(H_1(M))$ for any $M$ and that
Tor$(H_1(M))\neq T(H_1(M))$ for some $M$, e.g. 
$M=S^1\times S^1\times [0,1]$. 
Moreover, we note that, for a compact $3$-manifold $M$, 
$T(H_1(M))=H_1(M)$ if and only if $M$ can be embedded in 
a rational homology $3$-sphere. 
Let $K\cup K_1$ be a $2$-component oriented link in 
$M$ such that $K_1$ represents an element 
in Tor$(H_1(M))$, and let $c=|{\mathrm Tor}(H_1(M))|$. 
For a $2$-chain $F$ in 
$M$ with $\partial F=cK_1$, we define 
\[{\mathrm L}_M(K_1;K)\equiv{F\cdot K_2}\  
(\mbox{mod }{{\mathrm mul}([K])}).\]
Since $(F\cup (-F'))\cdot K$ is divisible by 
${\mathrm mul}([K])$ for any 2-chain $F'$ with 
$\partial F'=cK_1$, ${\mathrm L}_M(\ ;K)$ is well-defined. 
In the case that $[K]\in T(H_1(M))$, that is  
mul$([K])=0$, we may delete 
\lq modulo ${{\mathrm mul}([K])}$' from the definition 
above. If $[K]\in{\mathrm Tor}(H_1(M))$, then 
${\mathrm L}_M(K_1;K)/c=\lk_M(K_1,K)$.  $\Box$

\medskip
{\bf Remark 1.5.}
R.H. Kyle \cite{Kyl} showed that any symmetric integral matrix 
is congruent to a block sum of a nonsingular matrix 
and a zero matrix by integral unimodular matrix. 
This guarantees that any closed oriented $3$-manifold $M$ is obtained 
by Dehn surgery along a framed link 
$J_1\cup\cdots\cup J_n\cup J_1'\cup\cdots\cup J_m'$ 
in $S^3$ of which the linking matrix is a block sum 
of a nonsingular matrix $B$ and a zero matrix $O$,  
where $B$ (resp. $O$) is the linking matrix of $J_1\cup\cdots\cup J_n$
(resp. $J_1'\cup\cdots\cup J_m'$). 
By the arguments similar to that in proof of Theorem 1.1(2), we have  
the following: For a $2$-component oriented link $K_1\cup K_2$ in 
$S^3-J_1\cup\cdots\cup J_n\cup J_1'\cup\cdots\cup J_m'$ 
each of which component represents an element in Tor$H_1(M)$, 
\[\begin{array}{l}
\lk_M(K_1,K_2)-\lk_{S^3}(K_1,K_2)\\
\hspace*{3cm}=-(\lk_{S^3}(K_1,J_1),...,
\lk_{S^3}(K_1,J_n)){B}^{-1}
(\lk_{S^3}(K_2,J_1),...,\lk_{S^3}(K_2,J_n))^T. \Box
\end{array} \]

\bigskip
{\bf 2. Double branched cover of $S^3$}

\medskip
Let $K_1\cup K_1\cup\cdots\cup K_m$ be an $(m+1)$-component oriented link and 
$F$ and $F'$ unoriented surfaces bounded by $K$ 
without intersecting $K_1\cup\cdots\cup K_m$. 
These two surfaces are 
{\em $S^*$-equivalent rel. $K_1\cup\cdots\cup K_m$} if they 
are transposed into each other 
by the following operations; (1) attaching a half twisted band 
locally, (2) attaching a hollow 1-handle (1-surgery), and 
(3) deleting a hollow 1-handle (0-surgery),
where these operations can be done in the complement of 
$K_1\cup\cdots\cup K_m$. 

By the argument similar to that in the proof of Theorem 1 in \cite{Lee1}, 
we have the following.

\medskip
{\bf Proposition 2.1.} {\em Let $K\cup K_1\cup\cdots\cup K_m$ $(m\geq 1)$ be 
an oriented $(m+1)$-component link with 
the linking number $\lk(K,K_i)$ even for any $i(=1,...,m)$. 
Let $F$ and $F'$ be unoriented, possibly nonorientable surfaces bounded 
by $K$ without intersecting to $K_1\cup\cdots\cup K_m$. 
If $F$ and $F'$ are $S^*$-equivalent rel. $K_i\cup K_j$, then 
$\lambda_F(K_i,K_j)=\lambda_{F'}(K_i,K_j)$. }

\medskip
This theorem implies that $\lambda_F(K_i,K_j)$ is independent of the 
choice of a basis of $H_1(F)$.

\medskip
{\bf Remarks 2.2.} (1) Let $K\cup K_1\cup K_2$ be a split sum of 
a trivial knot $K$ and the Hopf link $K_1\cup K_2$. Let 
$F$ be a Seifert surface of $K$ as illustrated in 
Figure 1 and $D$ a disk bounded by $K$ with 
$D\cap(K_1\cup K_2)=\emptyset$. Then $\lambda_F(K_1,K_2)\neq 
\lambda_{D}(K_1,K_2)$. 
It follows from Proposition 2.1 that $F$ and $D$ are not 
$S^*$-equivalence rel. $K_1\cup K_2$. 
On the other hand, M. Saito \cite{Sai} showed that, 
for an oriented 2-component link $K\cup K_1$ with 
$\lk(K,K_1)$ even, any two unoriented surfaces bounded by $K$ without 
intersecting $K_1$ are $S^*$-equivalent rel. $K_1$. \\
(2)  In the next section, 
we will define $S$-equivalence, which is an orientable version of 
$S^*$-equivalence, rel. $K_1\cup K_2$ for Seifert surfaces for $K$ 
in $S^3\setminus(K_1\cup K_2)$. 
As a special case of \cite[Lemma 4]{M-M} or \cite[4.1.5 Proposition]{Kai}, 
we have that two Seifert surfaces for $K$ in $S^3\setminus(K_1\cup K_2)$ 
are $S$-equivalent rel $K_1\cup K_2$ if and only if they are homologous in 
$H_2(S^3\setminus(K_1\cup K_2), \partial N(K);{\Bbb Z})$, where $N(K)$ is a 
regular neighborhood of $K$ in $S^3\setminus(K_1\cup K_2)$. $\Box$

\begin{center}
\begin{tabular}{c} 
\includegraphics[trim=0mm 0mm 0mm 0mm, width=.25\linewidth]
{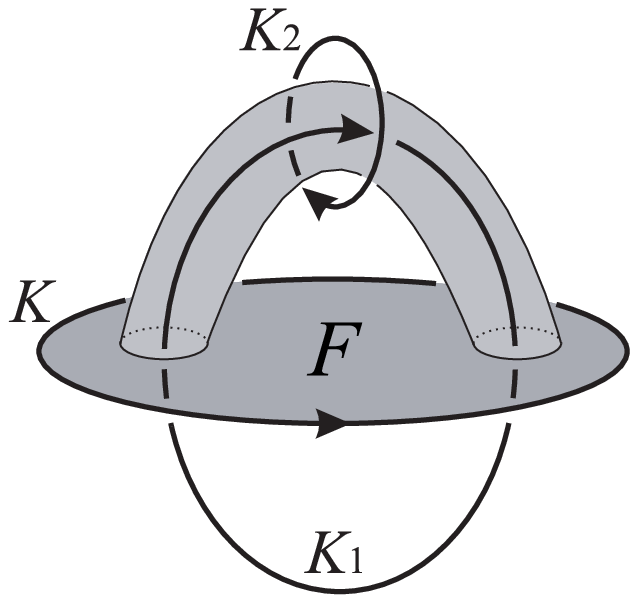}\\
Figure 1
\end{tabular}
\end{center}

\medskip
{\bf Proof of Proposition 2.1.}
Let $\beta$ be an another basis of $H_1(F)$. Then there is an unimodular 
matrix $P$ such that $\beta=\alpha P$, $G_{\beta}=P^T G_{\alpha}P$, 
$V_{\beta}(K_i)=V_{\alpha}(K_i)P$ and 
$V_{\beta}(K_j)^T=P^TV_{\alpha}(K_j)^T$. 
Thus we have 
\[V_{\beta}(K_i)G_{\beta}^{-1}V_{\beta}(K_j)^T=
V_{\alpha}(K_i)PP^{-1}G_{\alpha}^{-1}(P^T)^{-1}P^TV_{\alpha}(K_j)^T
=V_{\alpha}(K_i)G_{\alpha}^{-1}V_{\alpha}(K_j)^T.\] 

We may assume that 
$F'$ is obtained from $F$ by attaching a half twisted band 
or by attaching a hollow 1-handle. 

In the case that $F'$ is obtained from $F$ by attching a
half twisted band. Let $a$ be a cycle as illustrated in
Figure 2. Let $\alpha$ be a basis of $H_1(F)$ and 
$\beta=(a,\alpha)$ a basis of $H_1(F')$. 
Then we have 
\[G_{\beta}=\left(
\begin{array}{cccc}
 \pm1 & 0 & \cdots & 0\\
 0 &&& \\
 \vdots & & G_{\alpha} & \\
 0 &&&
\end{array}\right),\]
$V_{\beta}(K_i)=(0,V_{\alpha}(K_i))$ 
and $V_{\beta}(K_j)=(0,V_{\alpha}(K_j))$. 
Thus we have 
\[V_{\beta}(K_i)G_{\beta}^{-1}V_{\beta}(K_j)^T=
V_{\alpha}(K_i)G_{\alpha}^{-1}V_{\alpha}(K_j)^T.\]

Suppose that $F'$ is obtained from $F$ by attching a 
hollow 1-handle. Let $a$ and $b$ be cycles as illustrated in
Figure 3. Let $\alpha$ be a basis of $H_1(F)$ and 
$\beta=(a,b,\alpha)$ a basis $H_1(F')$. 
Then we have 
\[G_{\beta}=\left(
\begin{array}{ccccc}
0 & 1 & 0 & \cdots & 0\\
1 & x & x_1 & \cdots & x_n\\
0 & x_1 &&& \\
\vdots & \vdots & & G_{\alpha} & \\
0 & x_n &&&
\end{array}\right),\]
$V_{\beta}(K_i)=(0,\lk(K_i,b),V_{\alpha}(K_i))$ and 
$V_{\beta}(K_j)=(0,\lk(K_j,b),V_{\alpha}(K_j))$.
Then it is not hard to see that 
there are unimodular matrices $P$ and $Q$ such that 
\[P G_{\beta}Q=\left(
\begin{array}{ccccc}
0 & 1 & 0 & \cdots & 0\\
1 & 0 & 0 & \cdots & 0\\
0 & 0 &&& \\
\vdots & \vdots & & G_{\alpha} & \\
0 & 0 &&&
\end{array}\right),\]
$V_{\beta}(K_i)Q=V_{\beta}(K_i)$ and 
$P V_{\beta}(K_j)^T=V_{\beta}(K_j)^T.$
Thus we have 
\[V_{\beta}(K_i)G_{\beta}^{-1}V_{\beta}(K_j)^T=
V_{\beta}(K_i)(P G_{\beta}Q)^{-1}V_{\beta}(K_j)^T
=V_{\alpha}(K_i)G_{\alpha}^{-1}V_{\alpha}(K_j)^T.\]
This completes the proof. $\Box$

\begin{center}
\begin{tabular}{c}
\includegraphics[trim=0mm 0mm 0mm 0mm, width=.4\linewidth]
{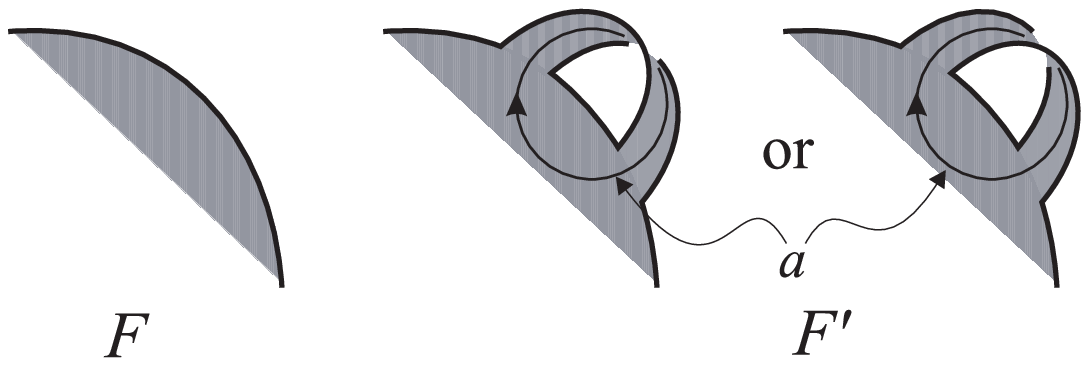}\\
Figure 2
\end{tabular}
\end{center}

\begin{center}
\begin{tabular}{c}
\includegraphics[trim=0mm 0mm 0mm 0mm, width=.6\linewidth]
{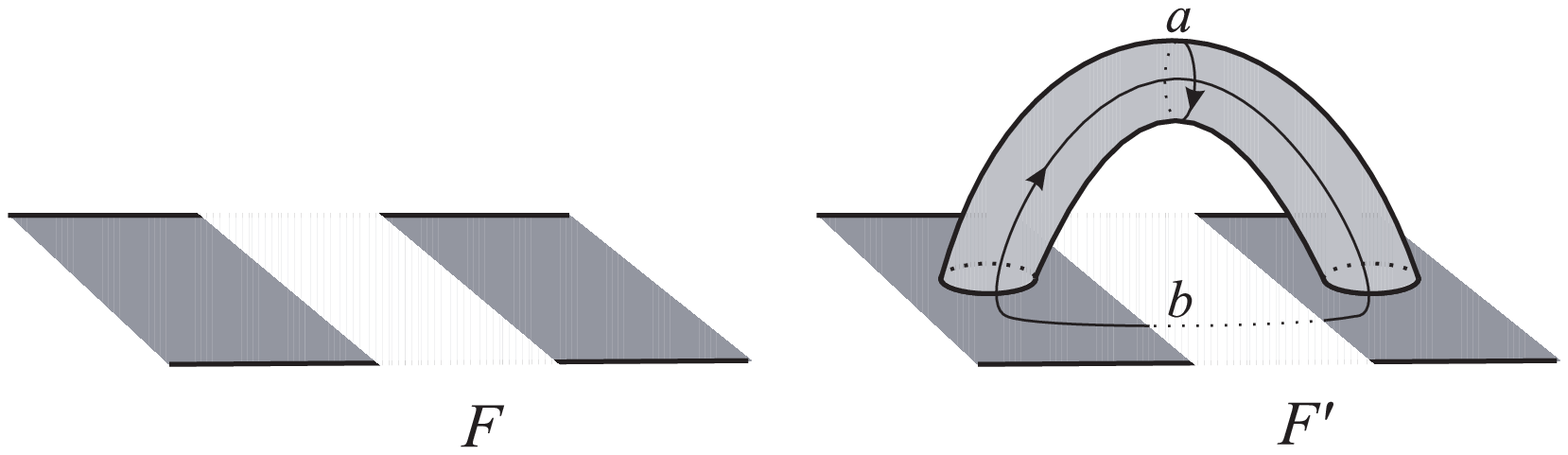}\\
Figure 3
\end{tabular}
\end{center}

\medskip
Let $K\cup K_1\cup K_2$ be an oriented $3$-component link 
with $\lk(K,K_i)$ even $(i=1,2)$. Let $F$ be an unoriented 
surface bounded by $K$ without intersecting $K_1\cup K_2$. 
According to the construction by S. Akbult and R. Kirby \cite{A-K}, 
we may assume that the double 
cover $X_2$ of $S^3$ branched over $K$ is obtained 
from two copies $M_1$ and $M_2$ of $S^3$-cut-along-$F$ by gluing 
each other along their boundaries suitably. 
Let $K_{1k}$ and $K_{2k}$ $(k=1,2)$ be the preimages in $M_k$ 
of $K_1$ and $K_2$, respectively. 

\medskip
{\bf Theorem 2.3.}  {\em For any $i,j,k,l$ $((i,k)\neq(j,l))$,
\[\lk_{X_2}(K_{ik},K_{jl})-(1-\delta_{ij})\delta_{kl}\lk(K_i,K_j)=
(-1)^{\delta_{kl}}\lambda_F(K_i,K_j),\]
where $\lk(K_i,K_i)=0$.}

\medskip
We note that  
$\lk_{X_2}(K_{i1},K_{i2})=\lambda_F(K_i,K_i)$ for each $i$ and 
$|\lk_{X_2}(K_{11},K_{21})-\lk_{X_2}(K_{11},K_{22})|=
|\lk_{X_2}(K_{12},K_{21})-\lk_{X_2}(K_{12},K_{22})|=
|2\lambda_F(K_1,K_2)-\lk(K_1,K_2)|$. 
Since $K_{ik}$'s are the preimage of $K_i$, we have the following 
corollary.  

\medskip
{\bf Corollary 2.4.} {Both $\lambda_F(K_i,K_i)$ $(i=1,2)$ and 
$|2\lambda_F(K_1,K_2)-\lk(K_1,K_2)|$ are  invariants of 
$K\cup K_1\cup K_2$. $\Box$}

\medskip
Now we denote $\lambda_F(K_i,K_i)$ by $\lambda_K(K_i)$.

\medskip
{\bf Remark 2.5.}
Let $K\cup K_1$ be an oriented link, $K_1(2,1)$ the $(2,1)$-cable knot 
of $K_1$. Since $\lk(K,K_1(2,1))$ is even, we can define 
\[\overline{\lambda}_K(K_1)=\frac{1}{4}\lambda_K(K_1(2,1)).\]
Note that $\overline{\lambda}_K(K_1)=\lambda_K(K_1)$ if 
$\lk(K,K_1)$ is even. Let $K\cup K_1\cup\cdots\cup K_m$ be an 
$(m+1)$-component oriented link. Then we define
\[\overline{\lambda}_K(K_1\cup\cdots\cup K_m)=\sum_{i=1}^m
\overline{\lambda}_K(K_i).\] 
Thus we have an invariant for oriented links. $\Box$

\medskip
{\bf Proof of Theorem 2.3.} 
Let $F$ be an unoriented surface bounded by $K$ with 
$F\cap(K_1\cup K_2)=\emptyset$. 
Then we may assume that $F$ is 
a surface as illustrated in Figure 4(a) or (b).
Let $a_i$ be a curve in $F$ as in Figure 4(a) or (b) $(i=1,...,n)$.  
Then we may regard that $(a_1,...,a_n)$ is a basis $\alpha$ of
$H_1(F)$. 
By \cite{A-K}, we have that the double branched cover $X_2$ 
is obtained from $S^3$ by 
Dehn surgery along an framed oriented link 
$J_1\cup\cdots\cup J_n$ with 
$\lk(J_i,J_j)=\lk(a_i,\tau a_j)$ for any $i\neq j$  
and with the framing of $J_i$ is equal to 
$\lk(a_i,\tau a_i)$ for any $i$. 
By the constraction, we note that
$K_{11}\cup K_{21}\cup K_{12}\cup K_{22}$ is in the 
complement of the framed link in $S^3$,  $\lk_{S^3}(K_{ik}, K_{jl})=
\delta_{kl}\lk_{S^3}(K_i,K_j)$ and 
\[(\lk_{S^3}(K_{ik},J_1),...,\lk_{S^3}(K_{ik},J_n))=
\left\{\begin{array}{ll}
V_{\alpha}(K_i) & \mbox{if $k=1$}\\
-V_{\alpha}(K_i) & \mbox{if $k=2$}.
\end{array}\right.\] 
By Corollary 1.2, we 
have the conclusion. $\Box$

\begin{center} 
\begin{tabular}{c}
\includegraphics[trim=0mm 0mm 0mm 0mm, width=.9\linewidth]
{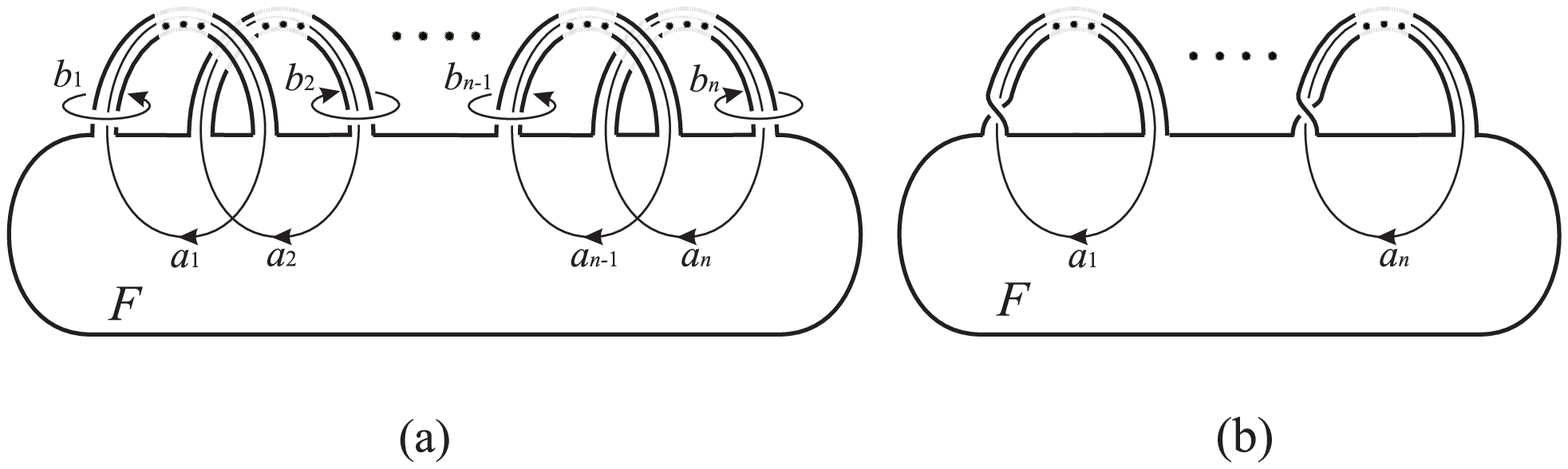}\\
Figure 4
\end{tabular}
\end{center}

\bigskip
{\bf 3. Cyclic branched cover of $S^3$}

\medskip
Let $K\cup K_1\cup\cdots\cup K_m$ be an $(m+1)$-component oriented link and 
$F$ and $F'$ Seifert surfaces of $K$ 
that donot intersect $K_1\cup\cdots\cup K_m$. 
These two surfaces are {\em $S$-equivalent rel. $K_1\cup\cdots\cup  K_m$} 
if they are transposed into each other 
by the following operations; (1) attaching a hollow orientable 1-handle 
(1-surgery), and (2) deleting a hollow 1-handle (0-surgery),
where these operations can be done in the complement of 
$K_1\cup\cdots\cup K_m$.

\medskip
By the argument similar to that in the proof of Proposition 2.1, 
we have the following theorem.

\medskip
{\bf Proposition 3.1.} {\em Let $K\cup K_1\cup\cdots\cup K_m$ $(m\geq 1)$ be 
an oriented $(m+1)$-component link with 
$\lk(K,K_i)=0$ for any $i(=1,...,m)$. 
Let $F$ and $F'$ be Seifert surfaces of $K$ 
that donot intersect to $K_1\cup\cdots\cup K_m$. 
If $F$ and $F'$ are $S$-equivalent rel. $K_i\cup K_j$, 
then $\lambda_F(K_i,K_j;\omega)=\lambda_{F'}(K_i,K_j;\omega)$, 
$\lambda_F(K_i,K_j)(t)=\lambda_{F'}(K_i,K_j)(t)$ and 
$\lambda_F^{(k,l)}(K_i,K_j)=\lambda_{F'}^{(k,l)}(K_i,K_j)$. $\Box$}

\medskip
This theorem implies that 
$\lambda_F(K_i,K_j;\omega)$,  
$\lambda_F(K_i,K_j)(t)$ and 
$\lambda_F^{(k,l)}(K_i,K_j)$ are independent of the 
choice of a basis of $H_1(F)$.

\medskip
Let $K\cup K_1\cup K_2$ be an oriented $3$-component link 
with $\lk(K,K_i)=0$ $(i=1,2)$. Let $F$ be a Seifert surface 
of $K$ with $F\cap(K_1\cup K_2)=\emptyset$. 
By \cite{A-K}, we may assume that the $p$-fold cyclic
cover $X_p$ of $S^3$ branched over $K$ is obtained 
from $p$ copies $M_1,...,M_p$ of $S^3$ by identifying 
$F\times[0,1]\subset M_i$ and $F\times[-1,0]\subset M_{i+1}$ suitably 
($i=1,...,p-1$), where $F\times\{0\}=F$. 
Let $K_{1k}$ and $K_{2k}$  be the preimages in $M_k$ $(k=1,...,p)$
of $K_1$ and $K_2$ respectively. 

\medskip
{\bf Theorem 3.2.} For any $i,j,k,l$ $((i,k)\neq(j,l))$, 
\[\begin{array}{l}
\lk_{X_{p}}(K_{ik},K_{jl})-(1-\delta_{ij})\delta_{kl}\lk(K_i,K_j)\\[2mm]
\hspace*{2.5cm}=\left\{
\begin{array}{ll}
-\lambda_F^{(k-1,l-1)}(K_i,K_j)+\lambda_F^{(k-1,l)}(K_i,K_j) &\\
\hspace*{2cm}+\lambda_F^{(k,l-1)}(K_i,K_j)-\lambda_F^{(k,l)}(K_i,K_j)
& \mbox{if $2\leq k\leq l \leq p-1$},\\[2mm]

\lambda_F^{(1,l-1)}(K_i,K_j)-\lambda_F^{(1,l)}(K_i,K_j) 
& \mbox{if $k=1,\ 2\leq l \leq p-1$},\\[2mm]

-\lambda_F^{(k-1,p-1)}(K_i,K_j)+\lambda_F^{(k,p-1)}(K_i,K_j)
& \mbox{if $2\leq k \leq p-1,\ l=p$},\\[2mm]

-\lambda_F^{(1,1)}(K_i,K_j)
& \mbox{if $k=l=1$},\\[2mm]

-\lambda_F^{(p-1,p-1)}(K_i,K_j)
& \mbox{if $k=l=p$},\\[2mm]

\lambda_F^{(1,p-1)}(K_i,K_j)
& \mbox{if $k=1$, $l=p$}.
\end{array}\right.
\end{array}\]

\medskip
{\bf Proof.} 
Let $F$ be a Seifert surface of $K$ with 
$F\cap(K_1\cup K_2)=\emptyset$. 
Then we may assume that $F$ is a surface as illustrated in Figure 4(a).
Let $a_i$ be a curve in $F$ as in Figure 4(a) $(i=1,...,n)$.  
Then we may regard that $(a_1,...,a_n)$ is a basis $\alpha$ of
$H_1(F)$. Let $M_{\alpha}$ be a Seifert matrix with respect to $\alpha$. 
By \cite{A-K}, we have that the $p$-fold cyclic branched cover $X_{p}$ 
is obtained from $S^3$ by 
Dehn surgery along an framed oriented link 
$J_{11}\cup\cdots\cup J_{n1}\cup\cdots\cup J_{1(p-1)}\cup
\cdots\cup J_{n(p-1)}$ 
with the linking matrix is equal to 
$M_{p,\alpha}$. 
By the constraction, we note that 
$K_{11}\cup K_{21}\cup\cdots\cup K_{ip}\cup K_{2p}$ is in the complement of 
the framed link in $S^3$, $\lk_{S^3}(K_{ik}, K_{jl})=
\delta_{kl}\lk_{S^3}(K_i,K_j)$ and 
\[\begin{array}{l}
(\lk_{S^3}(K_{ik},J_{11}),...,\lk_{S^3}(K_{ik},J_{n1}),...,
\lk_{S^3}(K_{ik},J_{1(p-1)}),...,\lk_{S^3}(K_{ik},J_{n(p-1)}))\\
\hspace*{7cm}=\left\{\begin{array}{ll}
V_{p,\alpha}^1(K_i) & \mbox{if $k=1$}\\
-V_{p,\alpha}^{p-1}(K_i) & \mbox{if $k=p$},\\
-V_{p,\alpha}^{k-1}(K_i)+V_{p,\alpha}^k(K_i) & \mbox{if $2\leq k\leq p-1$}.
\end{array}\right.
\end{array}\]  
By Corollary 1.2, we 
have the conclusion. $\Box$

\bigskip
{\bf 4. Infinite cyclic cover of the complement of a knot}

\medskip
{\bf Theorem 4.1.}
{\em Let $K\cup K_1\cup K_2$ be an oriented link with $\lk(K,K_1)=
\lk(K,K_2)=0$ and $F$ a Seifert surface of $K$ with $F\cap(K_1\cup K_2)=
\emptyset$. 
Let $(X_{\infty},p)$ be the infinite cyclic cover of the complement of $K$, 
$F_0$ a component of $p^{-1}(F)$, and  
${K_{ik}}$ a component of $p^{-1}(K_i)$ 
contained in a subspace bounded by 
$\tau^k F_0$ and $\tau^{k+1} F_0$ $(i=1,2,\ k\in {\Bbb Z})$.  
Then 
\[\tilde{\lk}_{X_{\infty}}({K_{1k}},{K_{2k}})
-\lk(K_1,K_2)=(1-t)\lambda_F(K_1,K_2)(t).\] 
Here $\tau$ is a covering translation that shifts $X_{\infty}$ 
along the positive direction of $p^{-1}(F)$. }

\medskip
Take a parallel copy $K'_i$ of $K_i$ in $S^3$ with $\lk(K_i,K'_i)=0$. 
Then we have $\tilde{\lk}_{X_{\infty}}({K_{ik}},{K'_{ik}})=
(1-t)\lambda_F(K_i,K'_i)(t)=(1-t)\lambda_F(K_i,K_i)(t)$.  
Meanwhile we note that
$\tilde{\lk}_{X_{\infty}}(\tau^m{K_{1k}},\tau^n{K_{2k}})=
t^{m-n}\tilde{\lk}_{X_{\infty}}({K_{1k}},{K_{2k}})
=t^{m-n}((1-t)\lambda_F(K_1,K_2)(t)+\lk(K_1,K_2)).$
Hence we have the following corollary.

\medskip
{\bf Corollary 4.2.}
{\em \begin{enumerate}
\item[$(1)$] $\lambda_F(K_i,K_i)(t)$ is an invariant of
$K\cup K_i$ and so is $\lambda_F(K_i,K_i;\omega)$.
\item[$(2)$] $(1-t)\lambda_F(K_1,K_2)(t)+\lk(K_1,K_2)$ 
is an invariant of $K\cup K_1\cup K_2$ 
up to multiplication by $t^{\pm n}$. $\Box$
\end{enumerate}
} 

\medskip
{\bf Remarks 4.3.} (1) As we mentioned in Introduction, 
for a parallel copy $K'_i$ of $K_i$ with $\lk(K_i,K'_i)=0$, 
the linking pairing $\tilde{\lk}_{X_{\infty}}(K_i,K'_i)$ is called  
Kojima-Yamazaki's $\eta$-function $\eta(K,K_i;t)$.
Thus $\eta(K,K_i;t)=(1-t)\lambda_F(K_i,K_i)(t)$, and hence
$\lambda_F(K_i,K_i)(t)$ is a topological concordance invariant. 
A different way to calculate the value of Kojima-Yamazaki's $\eta$-function 
was given in \cite{Jin}.  \\
(2) By the argument similar to that in 
\cite[Proof of Theorem 2]{K-Y}, we see that
$\tilde{\lk}_{X_{\infty}}({K_{1k}},{K_{2k}})(=
(1-t)\lambda_F(K_1,K_2)(t)+\lk(K_1,K_2))$ 
is a topological concordance invariant of $K\cup K_1\cup K_2$ 
up to multiplication by $t^{\pm n}$.\footnote{U. Kaiser pointed out that
the invariant $\tilde{\lk}_{X_{\infty}}({K_{1k}},{K_{2k}})$ was 
given by U. Dahlmeier \cite{Dah}.} 
%

\medskip
Let $\lambda_K(K_i)(t)$ and 
$\lambda_K(K_i;\omega)$ denote $\lambda_F(K_i,K_i)(t)$ and 
$\lambda_F(K_i,K_i;\omega)$ respectively. 
Note that 
$\lambda_K(K_i)(\overline{\omega})=(\omega-1)
\lambda_K(K_i;\omega)$.

\medskip
{\bf Remark 4.4.}
For an oriented 2-component link $K\cup K_1$ and 
for the untwisted double $K_1(2)$ of $K_1$, we define
\[\overline{\lambda}_K(K_1)(t)=\left\{
\begin{array}{ll} 
\lambda_K(K_1)(t) & \mbox{ if $\lk(K,K_1)=0$},\\
\lambda_K(K_1(2))(t) & \mbox{ otherwise}.
\end{array}\right.\]
For an $(m+1)$-component oriented link $K\cup K_1\cup\cdots\cup K_m$,  
we define
\[\overline{\lambda}_K(K_1\cup\cdots\cup K_m)(t)=\sum_{i=1}^m
\overline{\lambda}_K(K_i)(t).\]
Hence we have an invariant for oriented links. $\Box$

\medskip
By the definition of the linking pairing, we have the following 
lemma.

\medskip
{\bf Lemma 4.5.} {\em Let $X_{\infty}$ be the infinite cyclic cover 
of the complement of a knot and $K\cup K_1\cup\cdots\cup K_m$ $($resp. 
$K\cup K'_1\cup \cdots\cup K'_n)$ an oriented $(m+1)$-component 
$($ resp. $(n+1)$-component$)$ link in $X_{\infty}$. 
If there is a $2$-chain $F$ such that 
$\partial F=K_1\cup\cdots\cup K_m\cup(-K'_1)\cup\cdots\cup(-K'_n)$, 
then $\tilde{\lk}_{X_{\infty}}(K,K_1\cup\cdots\cup K_m)=
\tilde{\lk}_{X_{\infty}}(K,K'_1\cup\cdots\cup K'_n)+K\cdot F$ and 
$\tilde{\lk}_{X_{\infty}}(K_1\cup\cdots\cup K_m,K)=
\tilde{\lk}_{X_{\infty}}(K'_1\cup\cdots\cup K'_n,K)+K\cdot F$. 
Here $\tilde{\lk}_{X_{\infty}}(K,K_1\cup\cdots\cup K_m)=
\sum_{i=1}^m\tilde{\lk}_{X_{\infty}}(K,K_i)$ and 
$\tilde{\lk}_{X_{\infty}}(K_1\cup\cdots\cup K_m,K)=
\sum_{i=1}^m\tilde{\lk}_{X_{\infty}}(K_i,K)$. 
$\Box$}

\medskip
Let $K$ be a knot and $F$ a Seifert surface of $K$. 
We may assume that $F$ is a surface as illustrated in Figure 4(a). 
Let $a_1,...,a_n$ be curves as in Figure 4(a) and 
$M=(m_{ij})$ the Seifert matrix of $F$ with respect to 
a basis $[a_1],...,[a_n]$. 
Take curves $b_1,...,b_n$ so that $\lk(a_i,b_j)=\delta_{ij}$ 
for any $i,j$ as illustrated in Figure 4(a).  
Then we have the following lemma.

\medskip
{\bf Lemma 4.6. }{\em 
Let $(X_{\infty},p)$ be the infinite cyclic cover of the complement of $K$, 
$F_0$ a component of $p^{-1}(F)$, 
and ${b_{ik}}$ a component of $p^{-1}(b_i)$ 
contained in a subspace bounded by 
$\tau^k F_0$ and $\tau^{k+1} F_0$ $(i=1,...,n,\ k\in {\Bbb Z})$. 
Then $\lk_{X_{\infty}}({b_{ik}},{b_{jk}})$ 
is equal to the $(i,j)$-entry  of $(1-t)(tM-M^T)^{-1}$. }

\medskip
{\bf Proof.}
We denote by $t$ both a covering translation 
and a unit of ${\Bbb Z}[t,t^{-1}]$ since 
it is  well-known that there is natural correspondence between them.
Take curves $a^{\pm}_1,...,a^{\pm}_n$ so that 
$\lk(a_i^+,b_j)=\lk(a_i^-,b_j)=0$ for any $i,j$. 
Then $a_i^+$ is homologous to 
$\lk(a_i^+,a_1)b_1+\cdots+\lk(a_i^+,a_n)b_n$ and 
$a_i^-$ is homologous to 
$\lk(a_i^-,a_1)b_1+\cdots+\lk(a_i^-,a_n)b_n$. 
Moreover there are surfaces $E_i^+$ and $E_i^-$ 
that realize these homologous such that 
$E_i^+\cap F=E_i^-\cap F=\emptyset$ and 
$E_i^+$ (resp. $E_i^-$) is 
bounded by $-a_i^+$ (resp. $a_i^-$) and some copies of 
$b_j$'s $(j=1,...,n)$.
Then we have 
\[\left
(\begin{array}{c}
{[a_1^+]}\\
\vdots\\
{[a_n^+]}
\end{array}
\right)=M\left(
\begin{array}{c}
{[b_1]}\\
\vdots\\
{[b_n]}
\end{array}
\right),
\left
(\begin{array}{c}
{[a_1^-]}\\
\vdots\\
{[a_n^-]}
\end{array}
\right)=M^T\left(
\begin{array}{c}
{[b_1]}\\
\vdots\\
{[b_n]}
\end{array}
\right).\]
Let $A_i=a_i\times[-1,1]$ be an annulus in $S^3$ with 
$\partial A_i=\pm a_i\times\{\pm1\}=\pm a_i^{\pm}$ and ${A_{ik}}$  
a component of $p^{-1}(A)$ with $A_{ik}\cap t^{k+1}F_0
\neq \emptyset$. 
Then we have 
$\partial{A_{ik}}=t{a_{ik}^+}-{a_{ik}^-}$, 
where ${a_{ik}^{\pm}}$ is a component of $p^{-1}(a_i^{\pm})$ 
contained in a subspace between $t^kF_0$ and $t^{k+1}F_0$. 
Let ${E_{ik}^+}$ (resp. ${E_{ik}^-}$) 
be a component of $p^{-1}(E_i^+)$ (resp. $p^{-1}(E_i^-)$) 
contained in a subspace between $t^kF_0$ and $t^{k+1}F_0$. 
Let ${B_{ik}}={E_{ik}^-}\cup {A_{ik}}\cup{tE_{ik}^+}$. 
Then 
\[\left(
\begin{array}{c}
{[\partial{B_{1k}}]}\\
\vdots\\
{[\partial{B_{nk}}]}
\end{array}
\right)=
(t M-M^T)
\left(
\begin{array}{c}
{[{b_{1k}}]}\\
\vdots\\
{[{b_{nk}}]}
\end{array}
\right).\] 
Set $G(t)=tM-M^T$. 
Since $G(t)$ is nonsingular, we have
\[{\mathrm det}(G(t))G(t)^{-1}
\left(
\begin{array}{c}
{[\partial{B_{1k}}]}\\
\vdots\\
{[\partial{B_{nk}}]}
\end{array}
\right)=
{\mathrm det}(G(t))
\left(
\begin{array}{c}
{[{b_{1k}}]}\\
\vdots\\
{[{b_{nk}}]}
\end{array}
\right).\]
Set ${\mathrm det}(G(t))G(t)^{-1}=(l_{ij}(t))$. 
Since the boundary of each ${B_{ik}}$ is a disjoint union of 
some copies of ${b_{jk}}$'s and $t{b_{jk}}$'s $(j=1,...,n)$, 
$l_{i1}(t){B_{1k}}\cup\cdots\cup l_{in}(t){B_{nk}}$ is a 
$2$-chain of which boundary is a disjoint union of 
$t^s{b_{jk}}$'s ($s\in {\Bbb Z},\ j=1,...,n)$. 
Hence we have
\[\partial(l_{i1}(t){B_{1k}}\cup\cdots\cup l_{in}(t){B_{nk}})=
({\rm det}({G(t)}){b_{ik}})\cup
\bigcup_{1\leq j\leq n,s\in {\Bbb Z}}
c_{ijs}(t^s b_{jk}\cup(-t^s{b_{jk}})).\] 
Note that $\bigcup_{1\leq j\leq n,s\in {\Bbb Z}}
c_{ijs}(t^s{b_{jk}}\cup(-t^s{b_{jk}}))$ 
bounds a disjoint union $A$ of embedded annuli in
$X_{\infty}-p^{-1}(F)$. 
Since ${B_{ik}}\cdot {b_{jk}}=A_{ik}\cdot b_{jk}
=\delta_{ij}$, ${B_{ik}}\cdot t{b_{jk}}=A_{ik}\cdot t b_{jk}
=-\delta_{ij}$
 and ${B_{ik}}\cdot t^s{b_{jk}}=0$ 
for any $i,j$ and $s(\neq 0,1$), 
we have 
\[\lk_{X_{\infty}}({b_{ik}},{b_{jk}})=
\sum_{s\in{\Bbb Z}}\frac{t^s((A\cup 
l_{i1}(t){B_{1k}}\cup\cdots\cup l_{in}(t){B_{nk}})
\cdot t^s{b_{jk}})}{{\mathrm det}(G(t))}=
\frac{(1-t)l_{ij}(t)}{{\mathrm det}(G(t))}.\]
This completes the proof. $\Box$

\medskip
{\bf Proof of Theorem 4.1.}
It is not hard to see that there is a 2-component link 
$K_1\cup K_2'$ in $S^3-F$ such that
$K_2$ and $-K'_2$ cobound a surface $E_0$ in $S^3-F$ 
and $\lk(K_1\cup K_2')=0$. 
Let $E_{0k}$ (resp. $K_{2k}'$) be a component of $p^{-1}(E_0)$ 
(resp. $p^{-1}(K_2')$) contained in a subspace between 
$\tau^k F_0$ and $\tau^{k+1} F_0$.
Then, by Lemma 4.5, we have  
\[\begin{array}{rcl}
\tilde{\lk}_{X_{\infty}}(K_{1k},K_{2k})&=&
\tilde{\lk}_{X_{\infty}}(K_{1k},K_{2k}')+K_{1k}\cdot E_{0k}\\[2mm]
&=&\tilde{\lk}_{X_{\infty}}(K_{1k},K_{2k}')+K_1\cdot E_0\\[2mm]
&=&\tilde{\lk}_{X_{\infty}}(K_{1k},K_{2k}')+\lk(K_1,K_2).
\end{array}\]
Let $E_1$ (resp. $E_2$) be a Seifert surface of 
$K_1$ (resp. $K_2'$) in $S^3$ such that $E_1\cap K'_2=
K_1\cap E_2=\emptyset$ and 
$E_1\cup E_2$ intersects $F$ as illustrated in Figure 5. 
Let $N(F)$ be a small neighborhood of $F$ and 
$E^o_p=E_p-{\mathrm int}N(F)$ ($p=1,2$). 
Let $E^o_{pk}$ be a component of $p^{-1}(E^o_p)$ 
contained in a subspace between 
$\tau^k F_0$ and $\tau^{k+1} F_0$. 
We note that 
\[\partial E_{1k}^o=K_{1k}\cup\bigcup_{i=1}^n(-\lk(K_1,a_i)b_{ik}\cup
c_{i}b_{ik}\cup(-c_{i}b_{ik})),\]
and
\[\partial E_{2k}^o=K_{2k}'\cup\bigcup_{i=1}^n(-\lk(K'_2,a_i)b_{ik}\cup
d_ib_{ik}\cup(-d_ib_{ik})).\] 
By Lemma 4.5, we have 
\[\begin{array}{rcl}
\tilde{\lk}_{X_{\infty}}(K_{1k},K_{2k}')&=
&\displaystyle
\tilde{\lk}_{X_{\infty}}\left(\bigcup_{i=1}^n(\lk(K_1,a_i)b_{ik}\cup
(-c_ib_{ik})\cup c_ib_{ik}),K_{2k}'\right)\\[2mm]
&=&\displaystyle
\sum_{i=1}^n\lk(K_1,a_i)\tilde{\lk}_{X_{\infty}}(b_{ik},K_{2k}')\\[2mm]
&=&\displaystyle
\sum_{i=1}^n\lk(K_1,a_i)
\tilde{\lk}_{X_{\infty}}\left(b_{ik},\bigcup_{j=1}^n(\lk(K'_2,a_j)b_{jk}\cup
(-d_ib_{jk})\cup d_ib_{jk})\right)\\[2mm]
&=&\displaystyle
\sum_{i=1}^n\lk(K_1,a_i)\sum_{j=1}^n\lk(K'_2,a_j)
\tilde{\lk}_{X_{\infty}}(b_{ik},b_{jk}).\\[2mm]
\end{array}\]
Combining this and Lemma 4.6, we have 
\[\tilde{\lk}_{X_{\infty}}(K_{1k},K'_{2k})=(1-t)\lambda_F(K_1,K_2')(t).\]
Since $\lambda_F(K_1,K_2)(t)=\lambda_F(K_1,K_2')(t)$ and 
$\tilde{\lk}_{X_{\infty}}(K_{1k},K_{2k})=
\tilde{\lk}_{X_{\infty}}(K_{1k},K_{2k}')+\lk(K_1,K_2)$, 
we have the required result. $\Box$

\begin{center}
\begin{tabular}{c} 
\includegraphics[trim=0mm 0mm 0mm 0mm, width=.25\linewidth]
{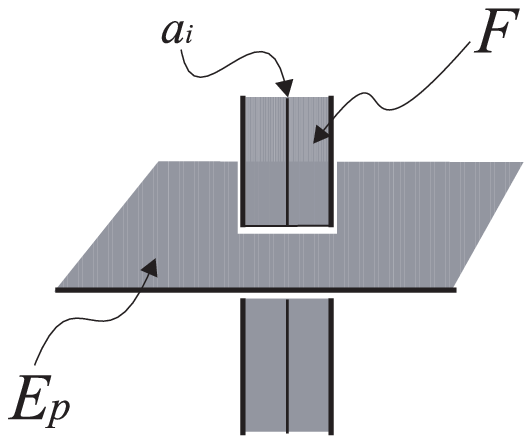}\\
Figure 5
\end{tabular}
\end{center}

\bigskip
{\bf 5. Connections between $\lambda_K$ and signatures}

\medskip
Let $K$ be a knot and $D$ a disk intersecting $K$ transversely in its 
interior with $|K\cap D|=2$. 
Performing $1/n$-Dehn surgery along $\partial D$, we obtain a new 
knot $K_n$. 
Note that if $\lk(\partial D,K)=0$ (resp. $\neq 0$) $K_{\pm1}$ 
(resp. $K_{\mp1}$) is obtained from $K$ by changing a $\mp$-crossing 
into a $\pm$-crossing. Then we have the following two theorems. 
These results were partially shown by Lee \cite{Lee1}, \cite{Lee2}. 
We modify his proofs. 

\medskip
{\bf Theorem 5.1.} {\em If $\lk(\partial D,K)=0$, then the following hold.
\begin{enumerate}
\item[$(1)$] 
$n(1-\omega)(1-\overline{\omega})\lambda_K(\partial D;\omega)\neq 1$ 
and $\lambda_K(\partial D,\omega)$ is a real number.
\item[$(2)$] $\sigma_{\omega}(K_n)=\sigma_{\omega}(K)-2n/|n|$ 
$($resp. $=\sigma_{\omega}(K))$ if and only if  
$n(1-\omega)(1-\overline{\omega})\lambda_K(\partial D;\omega)>1$ 
$($ resp. $<1)$. 
\item[$(3)$] $n(1-\omega)(1-\overline{\omega})\lambda_K(\partial D;\omega)= 
-\nabla_{K_n}(\sqrt{-1}|1-\omega|)/\nabla_{K}(\sqrt{-1}|1-\omega|)
+1$. 
\item[$(4)$] $n(t-1)\lambda_K(\partial D)(t)\neq 1$.
\item[$(5)$] 
$\Delta_{K_n}(t)$ is equal to 
$(1-n(t-1)\lambda_K(\partial D)(t))\Delta_K(t)$ 
up to multiplication by a unit of ${\Bbb Z}[t,t^{-1}]$.
\end{enumerate}
Here $\sigma_{\omega}$ is the Tristram-Levine signature 
{\rm\cite{Tri}}, $\nabla_K(z)$ is the Conway polynomial, and
$\Delta_K(t)$ is the Alexander polynomial. }

\medskip
By combining (2) and (3) in the theorem above, we have the following 
corollary.

\medskip
{\bf Corollary 5.2.} {\em Suppose $\lk(\partial D,K)=0$. 
Then $\sigma_{\omega}(K_n)=\sigma_{\omega}(K)-2n/|n|$ 
$($resp. $=\sigma_{\omega}(K))$ if and only if  
$\nabla_{K_n}(\sqrt{-1}|1-\omega|)/\nabla_{K}(\sqrt{-1}|1-\omega|)<0$ 
$($ resp. $>0)$. $\Box$}

\medskip
{\bf Remark 5.3.}
Note that a corossing change of a knot $K$ is realized by 
$\pm 1$-surgery along the boundary of disk $D$ with $|K\cap D|=2$ 
and $\lk(\partial D,K)=0$. 
By the corollary above and the induction on the unknotting number of a knot, 
we have that $\nabla_{K}(\sqrt{-1}|1-\omega|)/
|\nabla_{K}(\sqrt{-1}|1-\omega|)|=\sqrt{-1}^{\sigma_{\omega}(K)}$
for any knot $K$.  This implies that  
$\sigma_{\omega}(K)\neq 0$ if $\nabla_{K}(\sqrt{-1}|1-\omega|)<0$. 
$\Box$

\medskip
{\bf Theorem 5.4.} {\em 
\begin{enumerate}
\item[$(1)$] 
$2n\lambda_K(\partial D)\neq 1$.
\item[$(2)$] 
$\sigma(K_n)=\sigma(K)-2n/|n|+n|\lk(\partial D, K)|$ 
$($resp. $=\sigma(K)+n|\lk(\partial D, K)|)$ if and only if 
$2n\lambda_K(\partial D)>1$ 
$($ resp. $<1)$. 
\end{enumerate}
Here $\sigma(=\sigma_{-1})$ is the signature of a knot in the usual sense 
{\rm \cite{Tro1}, \cite{Mur}}
}

\medskip
{\bf Proof of Theorem 5.1.}
We note that there ia a Seifert surface $F$ of $K$ with 
$F\cap D$ is an arc as illustrated in Figure 6(a). 
We constract a Seifert surface $F'$ of $K_n$ from $F$ as 
illsutrated in Figure 6. 
Let $\alpha$ be a basis of $H_1(F)$. 
Let $a$ and $b$ be cycles as illustrated 
in Figure 6. 
We may assume that $\beta=(a,b,\alpha)$ is a basis of $H_1(F')$. 

Then we have  
\[G_{\beta,\omega}=\left(
\begin{array}{ccccc}
0 & 1-\overline{\omega} &  & 
\varepsilon(1-\overline{\omega})V_{\alpha}(\partial D) & \\
1-\omega & -n(1-\overline{\omega})(1-\omega) & 0  & \cdots & 0\\
 & 0 &&& \\
\varepsilon(1-\omega)V_{\alpha}(\partial D)^T & \vdots & & G_{\alpha,\omega} 
& \\
 & 0 &&&
\end{array}\right),\]
where $\varepsilon =1$ or $=-1$. This matrix is congruent to 
\[G_{\beta,\omega}'\left(
\begin{array}{ccccc}
1/n & 0 &  & 
\varepsilon(1-\overline{\omega})V_{\alpha}(\partial D) & \\
0 & -n(1-\overline{\omega})(1-\omega) & 0  & \cdots & 0\\
 & 0 &&& \\
\varepsilon(1-\omega)V_{\alpha}(\partial D)^T & \vdots & & G_{\alpha,\omega} & \\
 & 0 &&&
\end{array}\right).\]
Let  
\[U=\left(
\begin{array}{ccccc}
1 & 0 &  & 
-\varepsilon(1-\overline{\omega})V_{\alpha}(\partial D)
G_{\alpha,\omega}^{-1} & \\
0 & 1 & 0  & \cdots & 0\\
0 & 0 &&& \\
\vdots & \vdots & & I & \\
0 & 0 &&&
\end{array}\right).\]
Then we have 
\[UG_{\beta.\omega}^{-1}\overline{U}^T=
\left(\begin{array}{ccccc}
1/n-(1-\overline{\omega})(1-\omega)\lambda_K(\partial D;\omega) & 0 & 0 & 
\cdots & 0 \\
0 & -n(1-\overline{\omega})(1-\omega) & 0  & \cdots & 0\\
0 & 0 &&& \\
\vdots & \vdots & & G_{\alpha,\omega} & \\
0 & 0 &&&
\end{array}\right).\]
Thus $\lambda_K(\partial D,\omega)$ is a real number. 
Since this matrix is nonsingular, 
$\lambda_K(\partial D;\omega)\ne 1/n(1-\overline{\omega})(1-\omega)$. 
Moreover 
\[\sigma_{\omega}(K_n)={\mathrm sign}
\left(\begin{array}{cc}
1/n-(1-\overline{\omega})(1-\omega)\lambda_K(\partial D;\omega) & 0\\
0 & -n(1-\overline{\omega})(1-\omega) 
\end{array}\right)+\sigma_{\omega}(K).\]
This implies (1) and (2). 

Since (4) follows directly from (5), we shall prove (5). 
By the argument similar to that in the above, we have 
\[\begin{array}{rcl}
|G_{\beta}(t)|&=&\left|
\begin{array}{ccccc}
0 & t &  & 
\varepsilon tV_{\alpha}(\partial D) & \\
-1 & -n(t-1) & 0  & \cdots & 0\\
 & 0 &&& \\
-\varepsilon V_{\alpha}(\partial D)^T & \vdots & & G_{\alpha}(t) & \\
 & 0 &&&
\end{array}\right|\\
&&\\
&=& \left|
\begin{array}{ccccc}
-t/n(t-1)+t\lambda_K(\partial D)(t) & 0 & 0 & \cdots & 0 \\
0 & -n(t-1) & 0  & \cdots & 0\\
0 & 0 &&& \\
\vdots & \vdots & & G_{\alpha}(t) & \\
0 & 0 &&&
\end{array}\right|\\
&&\\
&=& t(1-n(t-1)\lambda_K(\partial D)(t))|G_{\alpha}(t)|.
\end{array}\]
Thus we have (5). 

In the proof of (5), replace $G_{\alpha}(t)$ and 
$G_{\beta}(t)$ with $t^{-1}M_{\alpha}-tM_{\alpha}^T$ and 
$t^{-1}M_{\beta}-tM_{\beta}^T$ respectively. 
By the argument similar to that in the proof of (5), we have 
\[\Omega_{K_n}(t)=(1-n(t^{-1}-t)
V_{\alpha}(\partial D)(t^{-1}M_{\alpha}-tM_{\alpha}^T)^{-1}
V_{\alpha}(\partial D)^T)\Omega_K(t),\]
where $\Omega_K(t)=|t^{-1}M_{\alpha}-tM_{\alpha}^T|$.
Put $t=\sqrt{-1}(1-\omega)/|1-\omega|$. Then we have 
\[\Omega_{K_n}\left(\frac{\sqrt{-1}(1-\omega)}{|1-\omega|}\right)
=(1-n(1-\omega)(1-\overline{\omega})\lambda_K(\partial D;\omega))
\Omega_{K}\left(\frac{\sqrt{-1}(1-\omega)}{|1-\omega|}\right).\]
Since $\Omega_K(t)=\nabla_K(t-t^{-1})$, we have 
\[\nabla_{K_n}(\sqrt{-1}|1-\omega|)
=(1-n(1-\omega)(1-\overline{\omega})\lambda_K(\partial D;\omega))
\nabla_{K}(\sqrt{-1}|1-\omega|).\]
The fact that $G_{\alpha,\omega}(=(1-\overline{\omega})M_{\alpha}+
(1-\omega)M_{\alpha}^T)$ is nonsingular implies  
$\nabla_{K}(\sqrt{-1}|1-\omega|)=
\Omega_{K}(\sqrt{-1}(1-\omega)/{|1-\omega|})
\neq 0$. Hence we have 
\[n(1-\omega)(1-\overline{\omega})\lambda_K(\partial D;\omega)= 
-\frac{\nabla_{K_n}(\sqrt{-1}|1-\omega|)}{\nabla_{K}(\sqrt{-1}|1-\omega|)}
+1.\]
This completes the proof.  $\Box$

\begin{center}
\begin{tabular}{c} 
\includegraphics[trim=0mm 0mm 0mm 0mm, width=.7\linewidth]
{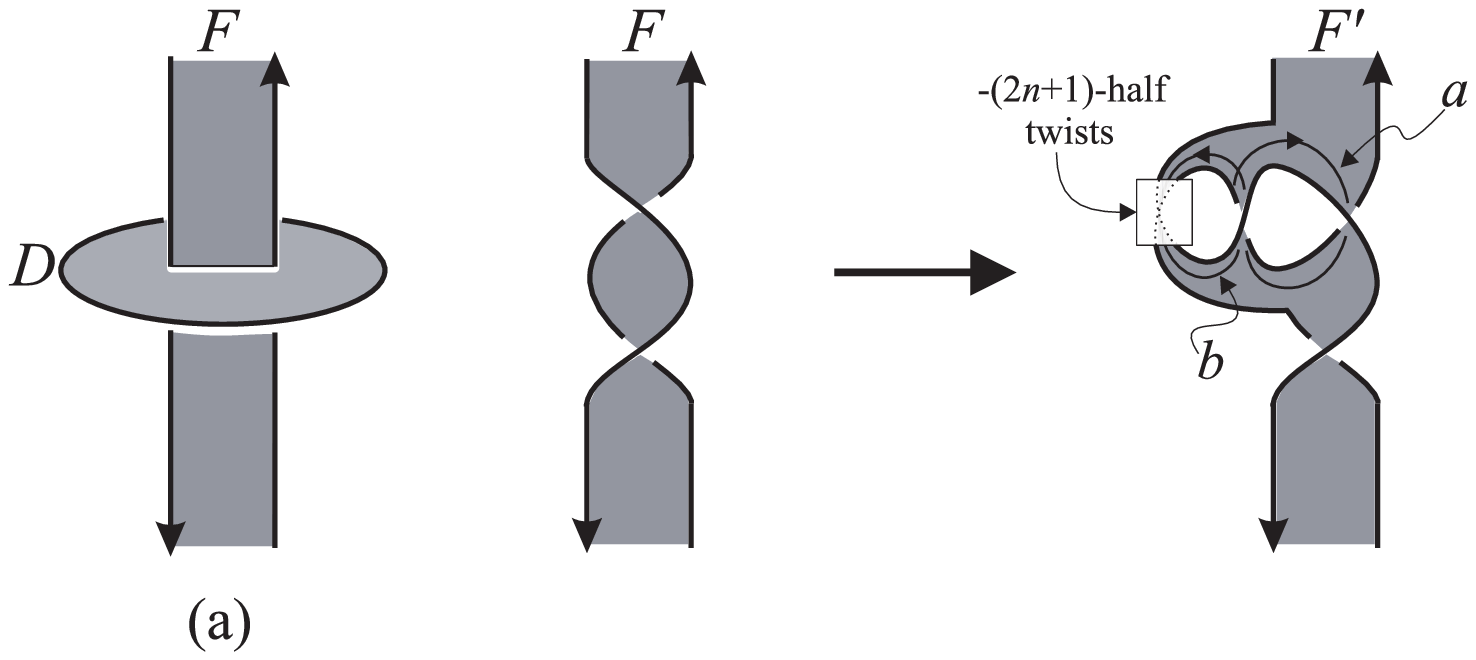}\\
Figure 6
\end{tabular}
\end{center}

\medskip
{\bf Proof of Theorem 5.4.}
Let $F$, $F'$, $\alpha$ and $\beta$ be the same as in the 
proof of Theorem 5.1. 
The only difference is that the surfaces are not necessarily orientable. 
Then we have  
\[G_{\beta,\omega}=\left(
\begin{array}{ccccc}
0 & 1 &  & 
\varepsilon V_{\alpha}(\partial D) & \\
1 & -2n & 0  & \cdots & 0\\
 & 0 &&& \\
\varepsilon V_{\alpha}(\partial D)^T & \vdots & & G_{\alpha} 
& \\
 & 0 &&&
\end{array}\right).\]
By the argument similar to that in the proof of Theorem 5.1, 
this matrix is congruent to 
\[
\left(\begin{array}{ccccc}
1/2n-\lambda_K(\partial D) & 0 & 0 & 
\cdots & 0 \\
0 & -2n & 0  & \cdots & 0\\
0 & 0 &&& \\
\vdots & \vdots & & G_{\alpha} & \\
0 & 0 &&&
\end{array}\right).\]
Since this matrix is nonsingular, 
$\lambda_K(\partial D)\neq 1/2n$. 
Moreover 
\[\sigma(K_n)={\mathrm sign}
\left(\begin{array}{cc}
1/2n-\lambda_K(\partial D) & 0\\
0 & -2n 
\end{array}\right)+{\mathrm sign}(G_{\alpha})+\frac{1}{2}e(F'),\]
where $e(F')$ is the {\em normal Euler number} of $F'$ \cite{G-L}. 
Since $e(F')=e(F)+2n|\lk(\partial D,K)|$,
\[{\mathrm sign}(G_{\alpha})+\frac{1}{2}e(F')=\sigma(K)+n|\lk(\partial D,K)|.\]
This completes the proof.
$\Box$


\bigskip
{\small
\begin{thebibliography}{9999}

\bibitem{A-K} S. Akbulut and R. Kirby, Branched covers of surfaces in 
$4$-manifolds, {\em Math. Ann.} {\bf 252} (1980), 111-131.

\bibitem{B-Z} G. Burde and H. Zieschang, {\em Knots}, 
De Gruyter Studies in Mathematics, {\bf 5},
Walter de Gruyter, Berlin, New York, 1985.


\bibitem{Dah} U. Dahlmeier, Verkettungshomotopien in Mannigfaltigkeiten, 
Doktorarbeit Siegen 1994.

\bibitem{Goe} L. Goeritz, Knoten und quadratische Formen, {\em Math. Z.} 
{\bf 36} (1933), 647-654.

\bibitem{G-L} C.McA. Gordon and R.A. Litherland, On the signature of a link, 
{\em Invent. Math.} {\bf 47} (1978), 53-69.

\bibitem{Hos} J. Hoste, A formula for Casson's invariant, 
{\em Trans. Amer. Math. Soc.} {\bf 297} (1986), 547-562.

\bibitem{Jin} G.T. Jin, On Kojima's $\eta$-function of links, 
{\em Differential topology}, pp. 14-30, Lecture Notes in Math. {\bf 1350}, 
Springer-Verlag, 1988.

\bibitem{Kai} U. Kaiser, {\em Link theory in manifolds}, 
Lecture Notes in Math. {\bf 1669}, Springer-Verlag, 1997.

\bibitem{Kau} L.H. Kauffman, Branched covering, open books and 
knot periodicity, {\em Topology} {\bf 13} (1974), 143-160.

\bibitem{Kea} C. Kearton, Blanchfield duality and simple knots, 
{\em Trans. Amer. Math. Soc.} {\bf 202} (1975), 141-160.

\bibitem{K-Y} S. Kojima and  M. Yamasaki, 
Some new invariants of links, 
{\em Invent. Math.} {\bf 54} (1979), 213-228.

\bibitem{Kyl} R.H. Kyle, 
Branched covering spaces and the quadratic forms of links, 
{\em Ann. of Math.} {\bf 59} (1954), 539-548.

\bibitem{Lee1} Y.W. Lee, A rational invariant for knot crossings, 
{\em Proc. Amer. Math. Soc.} {\bf 126} (1998), 3385-3392.

\bibitem{Lee2} Y.W. Lee, Alexander polynomial for link crossings, 
{\em Bull. Korean Math. Soc.} {\bf 35} (1998), 235-258.

\bibitem{Lev} J. Levine, Kont modules I, 
{\em Trans. Amer. Math. Soc.} {\bf 229} (1977), 1-50.

\bibitem{M-M} R. Mandelbaum and B. Moishezon, Numeric invariants in 
$3$-manifolds, {\em Low dimensional topology}, pp 285-304, 
Contemporary Math. {\bf 20} ed. S.J. Lomonaco jr., 
American Math. Soc. 1983. 


\bibitem{Mur} K. Murasugi, On a certain numerical invariant of link type, 
{\em Trans. Amer. Math. Soc.} {\bf 117} (1965), 387-422.


\bibitem{Rol} D. Rolfsen, {\em Knots and links}, 
Publish or Perish, Inc., Berkeley, 1976. 

\bibitem{Sai} M. Saito, On the unoriented Sato-Levine invariant, 
{\em J. Knot Theory Ramifications} {\bf 2} (1993), 335-358.

\bibitem{S-T} H. Seifert and W. Threlfall, 
{\em Lehrbuch der Topologie}, Teubner, Leipzig 1934.

\bibitem{Sok} M.V. Sokolov, Quantum invariants, skein modules, and 
periodicity of $3$-manifolds, Ph.D. Thesis, The George Washington 
University, Washington, D.C., 2000.

\bibitem{Tri} A. G. Tristram, Some cobordism invariants for links, 
{\em Proc. Cambridge Philos. Soc.} {\bf 66} (1969), 251-264.   	

\bibitem{Tro1} H.F. Trotter, Homology of group systems with applications 
to knot theory, {\em Ann. of Math.} {\bf 76} (1962), 464-498

\bibitem{Tro2} H.F. Trotter, On $S$-equivalence of Seifert matrices, 
{\em Invent. Math.} {\bf 20} (1973), 173-207.

\end{thebibliography} }
\end{document}